# A counterexample to linear relationship between largest common subtrees and smallest common supertrees[*]

Maciej Brzeski

Institute of Computer Science and Mathematics
Jagiellonian University
Kraków, Poland

Informatica

`maciej.brzeski@doctoral.uj.edu.pl`

We present a counter-example to the statement that there is a linear reduction in both directions between the largest common subtree and the smallest common supertrees. Moreover, we show that in general, even having of one of them, it is necessery to calculate the other from the linear number of nodes.

## 1  Introduction

Having two objects with tree-like structures, we often would like to know what is their common part, and what parts are distinct. When matching nodes, often we want to take into account the node neighborhood and tree structure. It can be defined in various ways [1], for example, by looking for the largest matching with certain conditions. Another approach is to define a tree editing operation, and find the shortest path to convert one tree into another [5], and this problem is called the tree edit distance. If we add the condition that we first add vertices and then remove them, we get the problem of tree alignment distance [3]. In this paper we discuss yet another approach, namely the largest common subtree and smallest common supertree [2].

However, it turns out that these definitions are closely related, and under certain assumptions, some of them become equivalent. In the tree edit distance problem, if we allow two operations: node insertion and node deletion, with the same cost for all nodes, is equivalent to the problem of finding the largest common subtree. It can be reasoned in this way, that the smallest edit distance we always obtain by first removing nodes from first tree, and then adding nodes from second tree. The tree which is after the last removal, and before the first insertion, is the largest common subtree.

In the paper, we will show that the construction of the smallest common supertree from the largest common subtree presented in [4] contains an error, and a counter-example showing that the proven relation between the two is false. The structure of this paper is as follows. In section 2 show a statement of the problem and a counter-example to it. In section 3 we extend a counter-example, showing that even having one of the largest common subtree or the smallest common supertree, does not reduce the complexity of computing the other. Section 4 summarizes results.

## 2  A counterexample

### 2.1  Definitions

A *directed graph* is a structure $G = (V, E)$, where $V$ is a set of nodes, and $E$ is a set of ordered pairs $(a, b) \in V \times V$ with $a \neq b$; the elements of $E$ are called arcs. For every arc $(v, w) \in E$, $v$ is its source node and $w$ its target node.

---

[*]This work is the result of research project no. DWD/4/66/2020 supported and funded by Ministry of Education and Science in Poland.

The in-degree of a node v in a graph is the number of arcs that have v as target node and its out-degree is the number of arcs that have v as source node.

A *path* in a directed graph $G = (V, E)$ is a sequence of nodes $(v_0, v_1, \ldots, v_k)$ such that $(v_0, v_1)$, $(v_1, v_2)$, $(v_2, v_3)$, ..., $(v_{k-1}, v_k) \in E$ ; its origin is $v_0$, its end is $v_k$, and its intermediate nodes are $v_1, \ldots, v_{k-1}$. Such a path is called non-trivial if $k \leq 1$. We shall represent a path from a to b, that is, a path with origin a and end b, by $a \rightsquigarrow b$.

A (*rooted*) *tree* is a directed finite graph $T = (V, E)$ with V either empty or containing a distinguished node $r \in V$, called the root, such that for every other node $v \in V$ there exists one, and only one, path $r \rightsquigarrow v$. Note that every node in a tree has in-degree 1, except the root that has in-degree 0. Henceforth, and unless otherwise stated, given a tree T we shall denote its set of nodes by $V(T)$ and its set of arcs by $E(T)$. The size of a tree T is its number $|V(T)|$ of nodes.

Now we can define what minor tree and minor embedding is, and then what is the largest common minor and smallest common supertree under minor embeddings.

**Definition 1.** *Let* S *and* T *be trees.* S *is a minor of* T *if there exists an injective mapping* $f: V(S) \to V(T)$ *satisfying the following condition: for every* $a, b \in V(S)$, *if* $(a, b) \in E(S)$, *then there exists a path* $f(a) \rightsquigarrow f(b)$ *in* T *with no intermediate node in* $f(V(S))$. *The mapping* f *is said to be a minor embedding* $f: S \to T$.

**Definition 2.** *Let S and T be trees.*

(i) *A largest common minor of S and T is a tree that is an minor of both of them and has the largest size among all trees with this property.*

(ii) *A smallest common supertree under minor embeddings of* S *and* T *is a tree such that both* S *and* T *are minor trees of it and has the least size among all trees with this property.*

In the rest of the paper, we will assume that all definitions apply to minor embeddings and will omit this for simplicity.

The construction in paper [4] is based on the quotient supergraph $T_{po}$, which we will define below, along with the mappings $\ell_1$ and $\ell_2$.

Let $T_1$ and $T_2$ be two trees. Let $T_\mu$ be a largest common subtree of them with $g_1: T_\mu \to T_1$ and $g_2: T_\mu \to T_2$ embeddings. Let $T_1 + T_2$ be the graph obtained as the disjoint sum of the trees $T_1$ and $T_2$: that is,

$$V(T_1 + T_2) = V(T_1) \sqcup V(T_2), \quad E(T_1 + T_2) = E(T_1) \sqcup E(T_2) \tag{1}$$

Let $\theta$ be the equivalence relation on $V(T_1) \sqcup V(T_2)$ defined, up to symmetry, by the following condition: $(a, b) \in \theta$ if and only if $a = b$ or there exists some $c \in V(T_\mu)$ such that $a = g_1(c)$ and $b = g_2(c)$.

**Definition 3.** *A* $T_{po}$ *is the quotient graph of* $T_1 + T_2$ *by equivalence* $\theta$ *with sets of nodes and arcs defined as follows:*

- $V(T_{po})$ *is the quotient set* $(V(T_1) \sqcup V(T_2))/\theta$, *with elements being the equivalence classes of the nodes of* $T_1$ *or* $T_2$;

- $E(T_{po})$ *is the set of arcs, induced by the arcs in* $T_1$ *or* $T_2$, *in the sense that* $([a], [b]) \in E(T_{po})$ *if and only if there exist* $a' \in [a]$, $b' \in [b]$ *and some i = 1, 2 such that* $(a', b') \in E(T_i)$.

Let $\ell_i: V(T_i) \to V(T_{po})$, $i = 1, 2$, denote the inclusion $V(T_i) \hookrightarrow V(T_1) \sqcup V(T_2)$ followed by the quotient mapping $V(T_1) \sqcup V(T_2) \to (V(T_1) \sqcup V(T_2))/\theta$: that is, $\ell_i(x) = [x]$ for every $x \in V(T_i)$. Note that, by construction,

$$V(T_{po}) = \ell_1(V(T_1)) \cup \ell_2(V(T_2)) \tag{2}$$

and

$$\ell_1(V(T_1)) \cap \ell_2(V(T_2)) = \ell_1(g_1(V(T_\mu))) = \ell_2(g_2(V(T_2))). \tag{3}$$

## 2.2 Problem statement (theorem from [4])

Let $T_1$ and $T_2$ be trees, and $T_\mu$ and $T_\sigma$ be the largest common subtree and smallest common supertree, respectively, along with $g_1\colon T_\mu \to T_1$, $g_2\colon T_\mu \to T_2$, $f_1\colon T_1 \to T_\sigma$ and $f_2\colon T_2 \to T_\sigma$ embeddings.

In the paper [4] an attempt has been made to prove that there exists a relation between largest common subtree and smallest common supertree, that allows each of them to be easily constructed from the other. In simple terms, the smallest common supertree can be obtained by merging nodes in sum of both trees, that are related to the same node of the largest common subtree. More formally, it was shown that the constructed supergraph $T_{po}$ in Def. 3 has certain properties that allow it to reduce into a supertree by removing all parallel arcs and all arcs subsumed by paths.

The error appears in the proof of the following proposition (numbering taken from the original work). In the original work, it was set for various types of embedding. Presented proposition is limited to minor embedding.

**Proposition 21.** ([4]) *Let $T_1$ and $T_2$ be trees, let $T_\mu$ be a largest common subtree of $T_1$ and $T_2$ with $g_1\colon T_\mu \to T_1$ and $g_2\colon T_\mu \to T_2$ embeddings, and let $T_{po}$ be the join of $T_1$ and $T_2$ obtained through $g_1$ and $g_2$.*

(i) *For every $v, w \in V(T_{po})$, if $(v,w) \in E(T_{po})$ and there is another path $v \rightsquigarrow w$ in $T_{po}$, then $v, w \in \ell_1(V(T_1)) \cap \ell_2(V(T_2))$, this path is unique and has no intermediate node in $\ell_1(V(T_1)) \cap \ell_2(V(T_2))$.*

(ii) *For every $v, w \in V(T_{po})$, if there are two different paths from $v$ to $w$ in $T_{po}$ without any common intermediate node, then one of them is the arc $(v, w)$, and then* (i) *applies.*

In particular, this construction implies that all nodes of the largest common subtree are also common in the smallest common supertree, as in Eq. 3. This also means that we have a equation between the sizes of these trees, which also appeared in the previous mentioned paper:

$$|T_\sigma| = |T_1| + |T_2| - |T_\mu|. \tag{4}$$

In [4] there were considered four types of embeddings. Our counterexample concerns minor embedding, as defined in Def. 1.

## 2.3 A counterexample

Consider trees $T_1$ and $T_2$ defined in Fig. 1a and 1c.

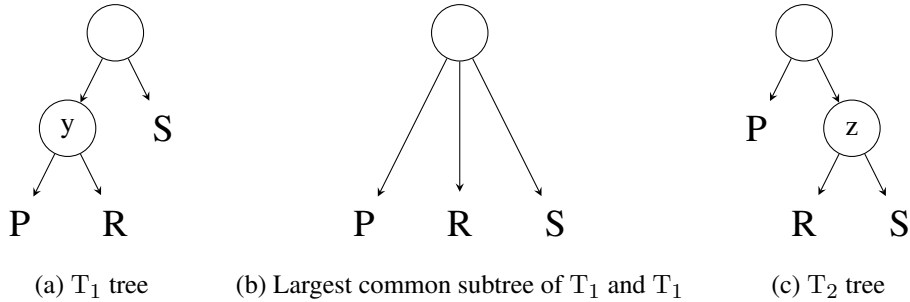

(a) $T_1$ tree     (b) Largest common subtree of $T_1$ and $T_1$     (c) $T_2$ tree

Figure 1: Largest common subtree

where P, R and S are any, non-empty trees. If P is not isomorfic with S, then the largest common subtree of $T_1$ and $T_2$ must have smaller number of nodes than $T_1$ (and $T_2$). The tree shown in Fig. 1b is the largest common subtree of $T_1$ and $T_2$.

In contrast, the smallest common supertree of $T_1$ and $T_2$ cannot be obtained in general by adding a single node to neither $T_1$ nor $T_2$. Now we will prove this, and this will contradict Eq. 4. First we prove below lemma:

**Lemma 4.** *Let S be a subtree of T with mapping $f\colon S \to T$. If $a, b \in V(S)$ and there is no path between $a$ and $b$, then there is also no path between $f(a)$ and $f(b)$.*

*Proof.* Suppose, without loss of generality, there is path from f(a) to f(b). Let c be lowest common ancestor of a and b in S. Note, than any path $x \rightsquigarrow y$, where $x, y \in S$, maps to a path $f(x) \rightsquigarrow f(y)$ in T. This occurs because we can split a path $x \rightsquigarrow y$ into arcs $(x, z_1), (z_1, z_2), \ldots, (z_n, y)$, all arcs map to paths $(f(x), f(z_1)), (f(z_1), f(z_2)), \ldots, (f(z_n), f(y))$, and finally we obtain a path concatenating all paths $f(x) \rightsquigarrow f(y) = ((f(x), f(z_1)), (f(z_1), f(z_2)), \ldots, (f(z_n), f(y)))$. In view of this, there must be a path $f(c) \rightsquigarrow f(a)$ and $f(c) \rightsquigarrow f(b)$ in T. However, the path $f(c) \rightsquigarrow f(b)$ must contain $f(a)$, because $f(c)$ in ancestor of $f(a)$, and there exists a path from $a \rightsquigarrow b$. We come to a contradiction, which ends the proof. □

**Theorem 5.** *Let $p_1$, $r_1$, $s_1$ be a nodes in trees $P_1$, $R_1$, $S_1$, respectively. Similarly, let $p_2$, $r_2$, $s_2$ be a nodes in trees $P_2$, $R_2$, $S_2$. Then it is impossible to all pairs, $p_1$ and $p_2$, $r_1$ and $r_2$, $s_1$ and $s_2$ be mapped into the same node of $T_\sigma$ simultaneously.*

*Proof.* Proof by contradiction. Let us assume that such a mapping is possible and $f_1(p_1) = f_2(p_2)$, $f_1(r_1) = f_2(r_2)$, $f_1(s_1) = f_2(s_2)$. Let us consider the relationship between nodes $f_1(y)$ and $f_2(z)$, shown in Fig. 1. There are three posibilities:

(i) If $f_1(y)$ is an ancestor of $f_2(z)$ (symmetrical case is analogous) then there would be a path between $f_1(y)$ and $f_1(s_1)$, which is impossible under lemma 4

(ii) If there is no path between $f_1(y)$ and $f_2(z)$ then both nodes would be the independent ancestors of $f_1(r_1) = f_2(r_2)$, which is impossible

(iii) If $f_1(y)$ and $f_2(z)$ are the same node in $T_\sigma$ then similarly there would be a path between $f_1(y)$ and $f_1(s_1)$, which is impossible under Lemma 4

None of the above cases can occur, which ends the proof. □

Using the above theorem, we can conclude with the following corollary:

**Corollary 6.** *In the example shown in Fig. 1 the smallest common subtree must be obtained by either adding one of the P, R or S trees in its entirety, or by creating the smallest common subtree of P and S trees, in two places*

All possibilities are presented in Fig. 2. The cases of adding P and S are symmetric, so both are included in case presented in Fig. 2a.

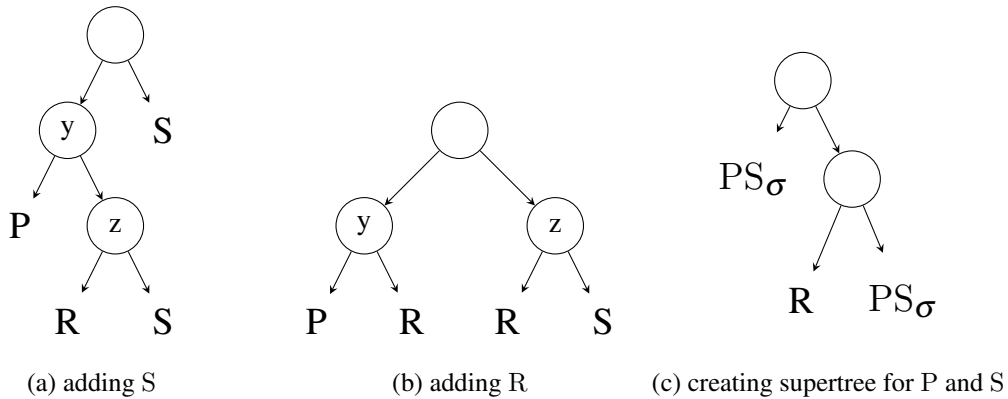

(a) adding S  (b) adding R  (c) creating supertree for P and S

Figure 2: Possible smallest common supertrees, where $PS_\sigma$ is the smallest common supertree for P and S

Note that in cases a) and b) the common supertrees are extended by one of the trees S, R or P and one extra node, which in total is more than one new vertex. For c), it is needed to expand both P and S to a common subtree. Since P and S can be any trees, so we can easily get the need to add more than one node.

In Proposition 21 from [4] the statement (ii) is incorrect. The reasoning behind the proof is that for trees defined in Fig. 1a and 1c, such a supergraph cannot contain a "diamond" shape, as in Fig. 3a. It is argued that then the largest common subtree can be extended by node x between a and R, and the mapping can by extended by $g_1(x) = y$ and $g_2(x) = y$, preserving the correctness of the embeddings. It has been shown that the created arcs $a \to x$ and $x \to \text{root}(R)$ map to paths without intermediate nodes in $g_1(T_\mu)$ or $g_2(T_\mu)$, but adding new mappings of $x \in T_\mu$ can cause that another path contains an intermediate node $g_i(x)$. In our example arc $a \to \text{root}(P)$ is mapped into path in tree $T_1$, where $y = g_1(x)$ is an intermediate node.

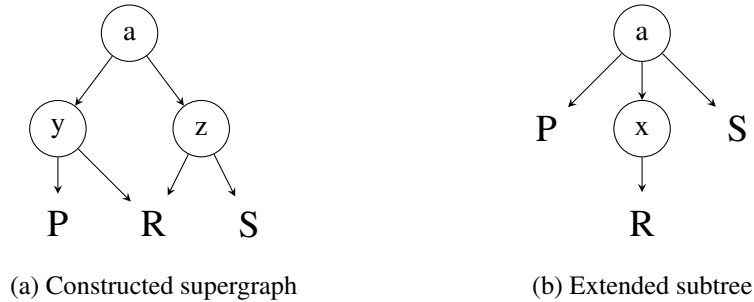

(a) Constructed supergraph    (b) Extended subtree

Figure 3: A counterexample for Proposition 21

## 3 Relationship between largest common subtree and smallest common supertree

In this section we will show examples, where even having one of the largest common subtree or the smallest common supertree, it is necessery to calculate the other from the linear number of nodes. Our examples are based on example given in Fig. 1.

Consider exactly such an example as in Fig. 1. Let P and S are defined as in Fig. 4, where A and B are any trees of sizes n and m, where without loss of generality $n \geq m$. R can be any tree with 2n nodes.

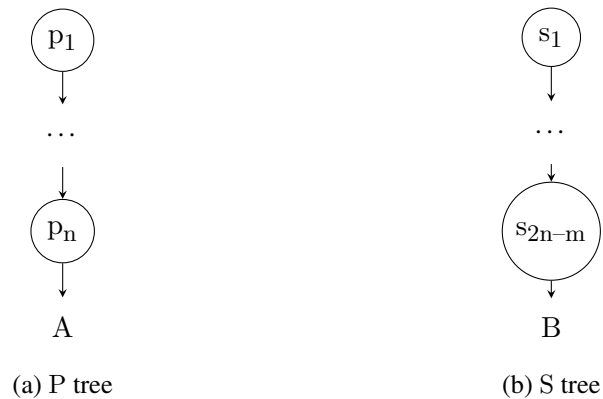

(a) P tree    (b) S tree

Figure 4

Considering the cases as in Fig. 1b we get that in case a) or b) we have to add $2n+1$ nodes (2n nodes from P, R or S and one node for not merging y and z), and at most 2n nodes in the case c), because all nodes $p_i$ and $s_i$ for $i = 1 \ldots n$ will be mapped into the same nodes in $T_\sigma$, and for the rest we need to add at most 2n nodes (the smallest common supertree cannot be larger than the sum of the both trees). Thus, the case c) is the smallest common supertree, where we have to calculate largest common supertree for A and B trees, which we chose arbitrarily.

We show now, that the same occurs in the opposite direction. Consider slightly changed $T_1$ and $T_2$ trees, defined in Fig. 5a and 5b. P and S trees are defined in Fig. 5c and 5d, and A and B are any trees of sizes n and m, where without loss of generality $n \geq m$.

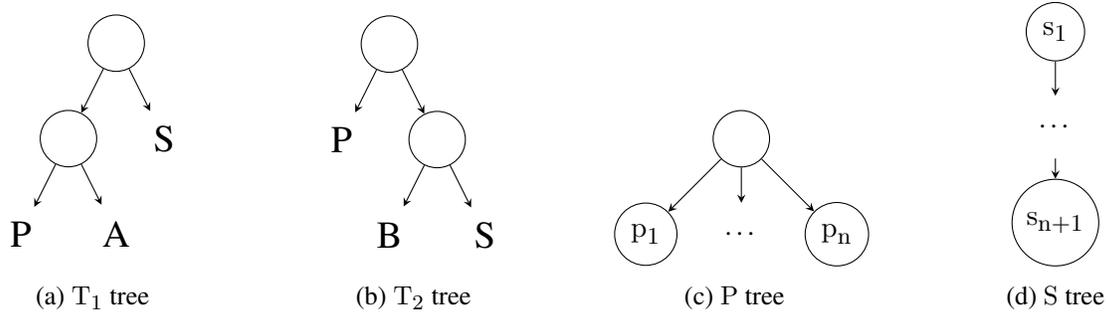

(a) $T_1$ tree  (b) $T_2$ tree  (c) P tree  (d) S tree

Figure 5

The smallest common supertree of P and S contains $2n-1$ nodes, because only two pairs of them can be mapped into the same node. Adding B increases the common supertree size by $m \leq n$, in other cases $n+1$ (adding S or P) or $2n-1$ nodes. Moreover, it is also needed to enlarge A to the common subtree of A and B. Thus, the smallest tree is obtained by adding B into $T_1$ (or A into $T_2$), and is presented in Fig. 6a. However, the largest common subtree has to be as in Fig. 6b. Finally, we have to calculate the largest common subtree for A and B, which we chose arbitrarily.

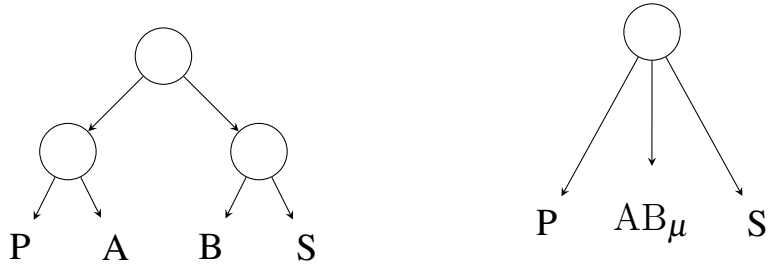

(a) Smallest common supertree for $T_1$ and $T_2$    (b) Largest common subtree for $T_1$ and $T_2$

Figure 6: Smallest common supertree and largest common subtree of $T_1$ and $T_2$, where $AB_\mu$ is the is the largest common subtree for A and B

## 4 Conclusions

We showed that, in general, the problem of the largest common subtree and the minimum common supertree, despite the apparent similarity of definitions, lead to completely different results. It is not possible to get one of them significantly faster with the result of the other. The differences are mainly due to the fact that the problem of the smallest common subtree is much more restrictive in terms of preserving the tree structure.

The above reasoning is also correct for labeled trees (it can be assumed that all vertices of an unlabeled tree have the same label). Since both problems are Max-SNP hard ([6, 3]), we can conclude that there is no polynomial reduction from either one into the other.